\documentclass[10pt]{article}
\usepackage{amsthm,amsmath,amsfonts}
\usepackage{fullpage,enumerate}
\usepackage{tikz}
\usepackage[colorlinks,citecolor=blue,urlcolor=blue]{hyperref}
\usepackage{hypernat}
\usepackage[notref,notcite]{showkeys}
\usepackage{graphicx}

\newcommand{\bZ}{\mathbb{Z}}

\newcommand{\la}{\lambda}

\begin{document}

\numberwithin{equation}{section}

\title{An inhomogeneous contact process model for speciation}

\bigskip
\author{  Rinaldo B. Schinazi
{\it University of Colorado, Colorado Springs}}
\maketitle

\bigbreak

{\bf Abstract.} We propose the following model on $\bZ^+$ for speciation and extinction. A species at site $i$ gives birth to a new species at site $j$ at rate $\la p(i,j)$ where $i$ and $j$ are nearest neighbors. A death at site $i$ occurs at rate $\delta_i$. We show that the existence of a phase transition in $\la$ depends critically on the value of the limit of $\frac{p(n,n+1)}{\delta_n}$.

\section{ The model}

Consider the following model on the positive integers $\bZ^+$.
Let $\lambda>0$ and $p$ be a transition probability matrix with the following properties.
$$ p(n,n+1)+p(n,n-1)=1\quad \forall n\geq 1,$$
and
$$ p(n,n+1)>0 \quad \forall n\geq 1.$$
For $n=0$ we assume
$p(0,1)=1.$
Let $(\delta_n)_{n\geq 0}$ be a sequence of strictly positive reals. We propose the following inhomogeneous contact process for speciation. We start the model with a single species at some $n\geq 0$.

(i) A species at $n\geq 1$ gives birth to another species at $n+1$ at rate $\lambda p(n,n+1)$ and gives  birth at $n-1$ at rate  $\lambda p(n,n-1)$. A species at $n=0$ can only give birth at 1 and this happens at rate $\lambda$.

(ii) A species at site $n\geq 0$ dies at rate $\delta_n>0$.

(iii) There is at most one species per site. Births on occupied sites are suppressed.

Moreover, we assume the existence of the following limit.

$$\lim_{n\to\infty} \frac{p(n,n+1)}{\delta_n}=\ell\leqno (1)$$
where $\ell\geq 0$ can possibly be infinite.

The basic contact process is the particular case $p(n,n+1)=p(n,n-1)=1/2$, $\delta_n=1$, for every $n$ in $\bZ$. See Liggett (1999) for a comprehensive mathematical survey of the contact process.

 As we will see below the probability that species survive forever depends critically on the value of $\ell$. The process is said to survive if there is a strictly positive probability that the process will survive. Otherwise it is said to die out.
\medskip

{\bf Theorem 1. }{\sl  Consider the inhomogeneous contact process on $\bZ^+$ starting with a single species.
\begin{itemize}
\item if $\ell=0$ then for any $\lambda>0$ the process dies out.
\item if $0<\ell<+\infty$ then there exists a critical value in $(0,+\infty)$ such that the process dies out for $\lambda$ below the critical value and survives for $\lambda$ above the critical value. 
\item if $\ell=+\infty$ then for any $\lambda>0$ the process survives.
\end{itemize}
}

\medskip

Biologically, an interesting particular case is when we have:

$$\lim_{n\to\infty} \delta_n=0\leqno (2)$$
and 
$$\lim_{n\to\infty} p(n,n+1)=0\leqno (3)$$

The number $1/\delta_n$ may be thought of as a measure of fitness for the species at $n$.
Hypothesis (2) assumes that as $n$ increases this measure of fitness goes to infinity. Hypothesis (3) assumes that it is harder and harder to give birth to a more fit species.  Hence, (2) and (3) go on opposite directions. Our result shows that extinction is certain for every $\la$ if and only if
$$p(n,n+1)<<\delta_n$$
for $n$ large. This relation can be interpreted as a cost-benefit trade-off.

Stochastic models for speciation go back to at least Yule (1924) who introduced a pure birth process for speciation. Birth and death processes are commonly used to model speciation, see Nee (2006) for a review of such models, see Aldous et al. (2011) for statistical questions related to birth and death models. Birth and death models are typically neutral: all the species are thought to have the same fitness. However, Liggett and Schinazi (2009) have introduced a birth and death chain with fitness. They assume that fitness is randomly assigned to a new species independently of its parent species, see also Guiol et al. (2011) and (2013). Independent fitnesses may be a reasonable assumption for virus but seems unreasonable for higher organisms. The present model can be seen as being at the other extreme of dependency since here the fitness of the new species depends heavily on the parent species.

\section{ Related models and results}

Inhomogeneous models for percolation and contact processes go back to at least Aizenmann and Grimmett (1991), see also Madras et al. (1994). The emphasis in those papers is to find conditions on the distribution of inhomogeneities that change the critical value of the process or make the phase transition discontinuous. There has also been significant work on the contact process in random environments, see Bramson et al. (1991), Liggett (1991) and (1992), Newman and Volchan (1996). The model introduced here is essentially the model from Liggett (1991). In that paper there are several sufficient conditions for extinction and a sufficient condition for survival.   However, Liggett's condition for survival is difficult to apply in general and seems to be useful only under very specific conditions such as periodic birth rates. Our approach is different in that we assume the existence of the limit $\ell$. This simplifies the problem and allows for a rather clear picture of survival and extinction for this particular case.

\section{ Proof of Theorem 1}

The proof is the consequence of two lemmas. The first one gives a sufficient condition for extinction.

\medskip

{\bf Lemma 1. }{ \sl If $\la\ell<1$ then the inhomogeneous contact process dies out.}

\medskip

We now prove Lemma 1. 
Let $\lambda>0$. We denote the inhomogeneous contact process by $\eta_t$.This process starts with a single species at some site $x$. Consider the following auxiliary process $\xi_t$ on $\bZ^+$. We start $\xi_t$ with a species at the same site $x$ and we fill all the sites of $\xi_0$ to the left of $x$ (if any). Let $r_t$ be the rightmost occupied site of $\xi_t$ at time $t$. The process $\xi_t$ evolves as follows. All the sites (if any) to the left of $r_t$ are always occupied. If $r_t=n$ then a death occurs at site $n$ at rate $\delta_n$ and a birth occurs at site $n+1$ at rate $\la p(n,n+1)$. 
These are the only transitions that are allowed for $\xi_t$. That is, a death can only occur at $r_t$ and a birth can only occur at $r_t+1$.
It is easy to see that the processes $\xi_t$ and $\eta_t$ can be coupled so that for any time $t\geq 0$ and any $x$ in $\bZ^+$ if $\eta_t(x)=1$ (i.e. $x$ is not empty) then $\xi_t(x)=1$. In particular, if $\xi_t$ dies out so does $\eta_t$. We now show that $\xi_t$ dies out.

 Observe that $r_t$ jumps from $n$ to $n+1$ at rate $\la p(n,n+1)$ and jumps from $n\geq 1$ to $n-1$ at rate  $\delta_n$. If $r_t=0$ the process dies out at rate $\delta_0>0$. Hence, $r_t$ is a birth and death chain on $\bZ^+$ with an absorbing state at 0.  A classical result (see for instance Karlin and Taylor (1975)) shows that $r_t$ is absorbed at 0 with probability 1 if and only if 
 $$\sum_{i=1}^\infty \frac{\delta_1\delta_2\dots\delta_i}{\la^ip(1,2)p(2,3)\dots p(i,i+1)}=+\infty.$$
 Since $\la\ell<1$ the ratio test implies that the series above is divergent. This completes the proof of Lemma 1.
 \medskip
 
 We now turn to a sufficient condition for survival. First we need to introduce the one-sided contact process. It is the particular case for which
 $$p(n,n+1)=1\mbox{ and }\delta_n=1\quad 
\forall n\geq 0.$$
There is a critical value $\la_c$ above which the one-sided contact process survives, see Harris (1976). 
 
 \medskip

{\bf Lemma 2. }{\sl Let $\la_c$ the critical value of the one-sided contact process on $\bZ$. If $\la\ell>\la_c$ then the inhomogeneous contact process on $\bZ^+$ survives.}

\medskip

We now prove Lemma 2.  Take $\la$ such that $\la\ell>\la_c$. For any  $\la'$ in $(\la_c,\la\ell)$ there is a positive integer $N$ such that if $n\geq N$ then
$$\la\frac{p(n,n+1)}{\delta_n}>\la'\leqno (4)$$

Consider now the one-sided contact process with birth rate $\la'$ and death rate 1 on the half-line $[N,+\infty)\cap \bZ$.  The one-sided contact process with birth rate $\la'$ and death rate 1 survives  on the half-line $[N,+\infty)\cap \bZ$. No births are allowed from $N-1$ to $N$ for this process. We now compare the inhomogeneous contact process to the one-sided contact process with the same boundary condition at site $N$. The rates being different the two processes run at different speeds. However, we can compare the two corresponding embedded discrete chains.  More precisely, we may compare the probabilities of a birth and the probability of a death. Using (4) we have that for $n\geq N$
$$\frac{\la p(n,n+1)}{\la p(n,n+1)+\delta_n}>\frac{\la'}{\la'+1}.$$
That is, a birth from $n$ to $n+1$ is more likely for the inhomogeneous contact chain than for the one-sided contact chain.  Note also that
$$\frac{\delta_n}{\la p(n,n+1)+\delta_n}<\frac{1}{\la'+1}.$$
That is, a death at $n$ is more likely for the one-sided contact chain than for the inhomogeneous contact process. Hence, births are more likely and deaths less likely for the inhomogeneous embedded chain. In particular survival of the one-sided process embedded chain implies survival of the inhomogeneous embedded chain.

Since the chain for the one-sided contact process has a positive probability of surviving so does the chain of the inhomogeneous contact process. Note that if $\ell=\infty$ then condition (4) holds for any $\la>0$. Hence, the process survives for any $\la>0$.
 This completes the proof of Lemma 2.
 
 \medskip
 
 We are now ready for the proof  of Theorem 1. 
 
 $\bullet$ If $\ell=0$ then $\la\ell<1$ for all $\la$ and Lemma 1 implies extinction.
 
  $\bullet$ If $0<\ell<\infty$ then Lemma 1 implies extinction for $\la\ell<1$ and Lemma 2 implies survival for $\la\ell>\la_c$. The existence of a critical value for the inhomogeneous contact process comes from the fact that the survival probability is increasing as a function of $\la$. This also implies that the critical value belongs to the interval $[1/\ell,\la_c/\ell]$.
  
   $\bullet$ If $\ell=\infty$  then Lemma 2 implies survival for every $\lambda$. This completes the proof of Theorem 1.
   \bigbreak
   
   {\bf References}
   
   D.J. Aldous, M.A. Krikun and L. Popovic (2011). Five statisitical questions about the Tree of Life. Systematic Biology Advance Access, published March 8 2011. 
   
   M. Aizenmann and G. Grimmett (1991). Strict monotonicity for critical points in percolation and ferromagnetic models. Journal of Statistical Physics 63, 817-835.
   
   M. Bramson, R. Durrett and R. Schonmann (1991). The contact process in a random environment. Annals of Probability 19, 960-983.
   
 H. Guiol, F. Machado and R.B. Schinazi (2011).  A stochastic model of evolution. Markov Processes Relat. Fields 17, 253-258.
 
  H. Guiol, F. Machado and R.B. Schinazi (2013). On a link between a species survival time in an evolution model and the Bessel distributions. Brazilian Journal of Probability and Statistics 
 27, 201-209. 
 
 T.E. Harris (1976). On a class of set-valued Markov processes. Annals of Probability 4, 175-194.
   
   S. Karlin and H. M. Taylor (1975). {\it A first course in stochastic processes.} Academic Press.
   
   T. Liggett (1991). Spatially inhomogeneous contact processes. In {\it Spatial stochastic processes. A Festschrift in honor of the seventieth birthday of Ted Harris}. K.C.Alexander and J.C. Watkins, editors. Birkhauser.
   
   T. Liggett (1992). The survival of one-dimensional contact processes in random environments. Annals of Probability 20, 696-723.
   
   T. Liggett (1999). {\it Stochastic interacting systems:contact, voter and exclusion processes.} Springer.
   
T.Liggett and R.B. Schinazi (2009).   A stochastic model for phylogenetic trees. Journal of Applied Probability 46 (2009), 601-607 . 
   
 N. Madras, R.B. Schinazi and R.H. Schonmann (1994). On the critical behavior of the contact process in deterministic inhomogeneous environments. Annals of Probability 22, 1140-1159.
   
   S. Nee (2006). Birth-Death models in macroevolution. Annu. Rev. Ecol. Evol. Syst. 37, 1-17.
   
   C.M.Newman and S. Volchan (1996). Persistent survival of one-dimensional contact processes in random environments. Annals of Probability 24, 411-421.
   
   G. Yule (1924). A mathematical theory of evolution based on the conclusions of Dr J.C. Willis. Philos. Trans. R. Soc. Lond. Ser. B. 213, 21-87.

\end{document}